\documentclass[12pt]{amsart}
\usepackage{amssymb, amsmath,amsthm,color,enumerate,booktabs,xcolor,cancel,framed,soul}
\usepackage{url,hyperref}
\usepackage[normalem]{ulem}
\hypersetup{citecolor=blue, linkcolor=blue, colorlinks=true}
\usepackage[marginparwidth=1.5cm,margin=2cm]{geometry}

\theoremstyle{plain}
\newtheorem{theorem}{Theorem}[section]
\newtheorem{lemma}[theorem]{Lemma}
\newtheorem{proposition}[theorem]{Proposition}
\newtheorem{corollary}[theorem]{Corollary}

\newtheorem{remark}[theorem]{Remark}

\newtheorem{defn}[theorem]{Definition}

\newcommand\PG{\mathsf{PG}}

\newcommand\s{\mathcal{S}}
\renewcommand\le{\leqslant}
\renewcommand\ge{\geqslant}

\definecolor{dkgreen}{rgb}{0,0.6,0}
\definecolor{gray}{rgb}{0.5,0.5,0.5}
\definecolor{mauve}{rgb}{0.58,0,0.82}
\definecolor{shadecolor}{rgb}{1,0.95,0.95}

\newcommand\q{\mathsf{Q}}
\newcommand\F{\mathbb{F}}
\newcommand\pg{\mathsf{PG}}
\newcommand\trace{\mathsf{Tr}}
\newcommand\norm{\mathsf{N}}
\newcommand\h{\mathsf{H}}
\newcommand\eq{\mathsf{Q}}
\newcommand\PGammaU{\mathrm{P\Gamma U}}

\title{On the Hermitian Veronesean}

\author{John Bamberg}
\address[Bamberg]{Department of Mathematics and Statistics, 
The University of Western Australia, 35 Stirling Highway, Perth, W. A. 6019, Australia.}

\author{Geertrui Van de Voorde}
\address[Van de Voorde]{School of Mathematics and Statistics, University of Canterbury, Private Bag 4800, 8140 Christchurch, New Zealand.}

\date{}



\begin{document}

\begin{abstract}
The Hermitian Veronesean in $\PG(3,q^2)$, given by  $\mathcal{V}:=\{ (1,x,x^q,x^{q+1}):x\in\F_q\}\cup\{(0,0,0,1)\}$, is a well-studied rational curve, and forms a {\em special} set of the Hermitian surface $\h(3,q^2)$. In this paper, we give two local characterisations of the Hermitian Veronesean, based on sublines and triples of points in perspective.
\end{abstract}

\maketitle

\section{Introduction}

Shult \cite{Shult:2005aa} introduced \emph{special sets} in connection with other objects in finite geometry,
such as generalised quadrangles, pseudo-ovals, and locally Hermitian ovoids.
A \emph{special set} of the Hermitian surface $\h(3,q^2)$, $q$ odd, is a set $\mathcal{S}$
of $q^2+1$ singular points such that for any point $P$ not in $\mathcal{S}$, there
are 0 or 2 points of $\mathcal{S}$ collinear with $P$. Equivalently, a special set of $\h(3,q^2)$ is a set of $q^2+1$ points such that any three span a nondegenerate plane.
The only known examples are equivalent under the collineation group of $\h(3,q^2)$  to the \emph{Hermitian Veronesean} $\mathcal{V}$ (described below), which arises as the $\F_{q^2}$-rational points of an 
$\F_{q^2}$-maximal degree $q+1$ rational curve (i.e., it has
genus 0 and meets the Hasse-Weil bound):
\[
\mathcal{V}:=\{ (1,x,x^q,x^{q+1}):x\in\F_q\}\cup\{(0,0,0,1)\}.
\]
For an excellent account of the properties of this curve, we refer to Lavrauw et al. \cite{Lavrauw:2023aa}. We will
refer to $\mathcal{V}$ as the \emph{classical} special set.
At the point of writing, it is still open whether or not every special set of $\h(3,q^2)$ is classical, and indeed, Hirschfeld and Thas list this problem as `Open problem 6' in their recent survey article \cite{Hirschfeld:aa}.

The main theorem of \cite{Cossidente:2006aa} (Theorem 3.1) states that a special set of $\h(3,q^2)$, such
that any three elements are \emph{in perspective} is classical. Here we use the term `in perspective'
for points of $\h(3,q^2)$ in accordance with its usual use in the dual setting (in $\q^-(5,q)$). While the result is correct, there is an error in the proof of their result whereby a frame of four points is chosen, when the hypothesis only allows three to be chosen. The same theorem was proved in a more general setting in 2011 by Thas \cite{Thas:2011aa},
where he showed that a pseudo-oval of $\pg(5,q)$ is elementary if and only if every triple is in perspective.
(Note that the dual of a special set is a pseudo-oval of $\pg(5,q)$ whose elements are lines of some 
elliptic quadric $\q^-(5,q)$.) A different proof for special sets was given in \cite[Corollary 5.3]{Bamberg:2021aa}. 

In this paper, we give two local characterisations of $\mathcal{V}$. The first one only requires that all triangles containing one of the three edges $PQ_i$, $i=1,2,3$
are in perspective: 

\begin{theorem}\label{main1}
Let $\s$ be a set of mutually noncollinear points on $\h(3,q^2)$, $q$ odd.
Suppose there are four non-coplanar points $P,Q_1,Q_2,Q_3\in \s$, such that all triangles $PQ_iR$, $i=1,2,3$, $R\in \s\setminus \{P,Q_i\}$, are in perspective.
Then $|\s|\le q^2+1$ and equality holds if and only if $\s$ is equivalent under the collineation group of $\h(3,q^2)$ to the Hermitian Veronesean $\mathcal{V}$.
\end{theorem}

Our second main theorem is a characterisation of $\mathcal{V}$ based on Baer sublines:

\begin{theorem}\label{main2}
    Let $q$ be an odd prime power, and let $\s$ be a set of $q^2+1$ pairwise noncollinear points of $\h(3,q^2)$ such that the following hold:
    \begin{enumerate}[(i)]
        \item There is a point $P\in\s$ and totally isotropic line $\ell$ through $P$, such that every point of $X$ incident
        with $\ell$ not equal to $P$ is collinear with a unique element of $\s\backslash\{P\}$. This yields
        a natural correspondence $F_\ell$ from $\ell$ to $\s$.
        \item Every Baer subline $b$ of $\ell$ containing $P$ maps under $F_\ell$ to $q+1$ 
        points of a plane of $\PG(3,q^2)$, incident with $P$.
        \item There is a point $Q$ such that for every $R\neq P,Q$ in $\s$, the triange $PQR$ is in perspective.
        \item The plane spanned by three distinct points $R_1,R_2,R_3\in \s$ is nondegenerate.
    \end{enumerate}
    Then $\s$ is equivalent under the collineation group of $\h(3,q^2)$ to the Hermitian Veronesean $\mathcal{V}$.   
\end{theorem}

As a consequence, we find that a special set satisfying conditions (ii) and (iii) is equivalent to the Hermitian Veronesean (see Corollary \ref{cor:sublinespecial}).



\section{Setup}
Throughout, $q$ is an odd prime power. We fix the following Hermitian form on a 4-dimensional vector space over $\F_{q^2}$:
\[
h(X,Y):=X_0Y_3^q+X_3Y_0^q-X_1Y_1^q-X_2Y_2^q,
\]
and let $\h(3,q^2)$ denote the associated Hermitian surface in $\PG(3,q^2).$
Shult \cite{Shult:2005aa} introduced the following function on triples of points of $\h(3,q^2)$, which
is closely related to the \emph{Segre-invariant} of three points $A$, $B$, $C$ (see \cite{Cossidente:2006aa}):
\[
[A,B,C]:=h(A,B)h(B,C)h(C,A)\F_{q}^*.
\]

\begin{remark}
    Segre's invariant on triples of points of $\h(3,q^2)$ is defined in \cite[\S53]{Segre65} to
    be the ratio 
    \[
    \frac{h(A,B)h(B,C)h(C,A)}{h(B,A)h(C,B)h(A,C)}.
    \]
    In particular, this quantity does not depend on the representative vectors of the points involved.
    We will adopt something similar to Shult's simpler invariant. Given $3$ points $A,B,C$ on $\h(3,q^2)$, we will slightly abuse notation and say that $[A,B,C]=h(A,B)h(B,C)h(C,A)$ is {\em the} Segre invariant of $A,B,C$. It is clear that this Segre invariant $[A,B,C]$ is only well-defined up to an $\F_q$-multiple. This is not an issue since we will only be concerned with Segre invariants being contained in $\F_q$ or having zero trace, both of which are invariant under $\F_q$-multiples. Similarly, we see that the Segre invariant of the points $A,B,C$ is equal to the Segre invariant of the points $B,A,C$ up to the action of the Frobenius automorphism $\psi:x\to x^q$. Since both $\F_q$ and $\{x\in \F_{q^2}\mid \trace(x)=0\}$ are stabilised by $\psi$, we can consider the triples of points under consideration as being unordered.
\end{remark}

Throughout, we will denote by $\trace$ the relative trace map from $\F_{q^2}$ to $\F_q$,
where $\trace(x)=x+x^q$ for all $x\in \F_{q^2}$. Likewise, $N(x):=x^{q+1}$ is the relative norm map
from $\F_{q^2}$ to $\F_q$.

\begin{lemma}[{\cite[Lemma 1]{Shult:2005aa}}, {\cite[Lemma 4.1]{Bamberg:2021aa}}]\samepage\label{lemma:degplane}
Let $P,Q,R$ be three distinct noncollinear points of $\h(3,q^2)$. Then
$P,Q,R$ spans a degenerate plane if and only if  
$\trace[P,Q,R]=0$.
\end{lemma}

Using $\perp$ to denote the image under the polarity associated with $\h(3,q^2)$, we say that three pairwise noncollinear points $A$, $B$, $C$ of $\h(3,q^2)$ are \emph{in perspective} if
the planes $\langle A^\perp\cap B^\perp, C\rangle$, $\langle B^\perp\cap C^\perp, A\rangle$, 
$\langle C^\perp\cap A^\perp, B\rangle$ intersect in a line. This is the dual of the more natural notion of perspective triples of lines on the elliptic quadric $\eq^-(5,q)$ (via the Klein
correspondence). We also say that the \emph{triangle} $ABC$ is in perspective.
 
\begin{lemma}[{\cite[Proposition 5.1]{Bamberg:2021aa}}]\label{lemma:inperspective}
Let $P,Q,R$ be three distinct noncollinear points of $\h(3,q^2)$. Then $PQR$ is in perspective if and only if $[P,Q,R]\in \F_q$.
\end{lemma}

As mentioned in the introduction,
Thas \cite{Thas:2011aa} showed that if there are $q^2+1$ lines $\mathcal{L}$ of $\PG(5,q)$, $q$ odd, such that any three are in perspective, then $\mathcal{L}$ is projectively equivalent to the
field reduction of a conic of $\PG(2,q^2)$. It is still not known if every \emph{pseudo-oval}
of $\PG(5,q)$ is the field reduction of a conic of $\PG(2,q^2)$ (a \emph{pseudo-conic}), and a pseudo-oval is a weaker object of $q^2+1$ lines in $\PG(5,q)$ having any three in general position.
From the interest of special sets, the pseudo-conic is the only example
of a set of $q^2+1$ lines of the elliptic quadric $\q^-(5,q)$ such that every line
of $\q^-(5,q)$ outside of it is concurrent with 0 or 2 elements within it.
We will be interested in the stronger form of Thas' Theorem, but for lines
of $\q^-(5,q)$. Below is the dual result.


\begin{theorem}[cf. {\cite[Theorem 6.4]{Thas:2011aa}}] \label{thm:perspective}
    Let $\s$ be a set of points on $\h(3,q^2)$ such that all triangles in $\s$ are in perspective. Then $\s$ is equivalent under the collineation group of $\h(3,q^2)$ to the Hermitian Veronesean.
\end{theorem}

\begin{corollary}\label{justonepoint}
    Let $\s$ be a set of points on $\h(3,q^2)$ such that there exists a point $P$ such that all triangles on $P$ are in perspective. Then $\s$ is equivalent under the collineation group of $\h(3,q^2)$ to the Hermitian Veronesean.
\end{corollary}

\begin{proof}
We use Lemma \ref{lemma:inperspective}.
Note that $[R_1,R_2,R_3]=[P,R_1,R_2]\cdot [P, R_2, R_3]\cdot [P, R_3, R_1]$, so if all the terms on the right-hand side lie in $\F_q$, then the left-hand side lies in $\F_q$.
\end{proof}

Note that Thas' result requires all triangles to be in perspective, whereas Corollary \ref{justonepoint} improves
this to `all triangles through a given point'. In Theorem \ref{main1}, we 
show that we need even fewer triangles -- we only demand that all triangles $PQ_iR$ are 
in perspective, for four non-coplanar fixed points $P, Q_1, Q_2, Q_3\in \s$. Moreover,
Theorem \ref{main1} does not depend on Corollary \ref{justonepoint}, and so it
reproves Theorem \ref{thm:perspective} and Corollary \ref{justonepoint} at the same time.

The following results give us the necessary symmetry results that allow us to
derive a standard form for the Hermitian Veronesean.

\begin{lemma}\label{lemma:stabPQS}
Let $\h(3,q^2)$, where is $q$ is a power of an odd prime $p$, be the Hermitian polar space with form defined by 
$h(X,Y):=X_0Y_3^q+X_3Y_0^q-X_1Y_1^q-X_2Y_2^q$. Let $P, Q, S$ be the points with homogeneous
coordinates $(0,0,0,1)$, $(1,0,0,0)$, $(1,1,1,1)$ respectively, of $\h(3,q^2)$.
Then the pointwise stabiliser $C$ in $\PGammaU(4,q)$ of $P,Q,S$ consists of elements whose
representative matrices are of the form
\[
\begin{bmatrix}
 1&0&0&0\\
 0&x&1-x&0\\
 0&1-x&x&0\\
 0&0&0&1
\end{bmatrix}\cdot \phi
\]
where $x+x^q=2x^{q+1}$ and $\phi$ is induced by a field automorphism of $\F_{q^2}$. 
Moreover, if $D$ is the setwise stabiliser in $\PGammaU(4,q)$ of $P,Q,S$,
then $D/C\cong S_3$, and the size of $D$ is $12(q+1)\log_p(q).$
\end{lemma}

\begin{proof}
First, let $g$ be a semilinear map $M \phi$, represented by an invertible matrix $M$
and field automorphism $\phi$, fixing each of the points $P$, $Q$, $S$ (with the action of $g$
on the right of these elements). Since $\phi$ fixes $P$, $Q$, $S$, it follows that
$M$ is, up to a scalar, of the form
\[
\begin{bmatrix}
    1 & 0 & 0 & 0\\
    a & b & c & d\\
    e & f & h & i\\
    0 & 0 & 0 & \mu
\end{bmatrix}
\]
such that 
\begin{equation}\label{eq:first}
    1+a+e=b+f=c+h=d+i+\mu.
\end{equation}
The Gram matrix of the Hermitian form
is
\[
U:=
\begin{bmatrix}
    0 & 0 & 0 & 1\\
    0 & -1 & 0 & 0\\
    0 & 0 & -1 & 0\\
    1 & 0 & 0 & 0
\end{bmatrix}.
\]

In order for $g$ to preserve the Hermitian form, there must exist a nonzero scalar $\lambda\in\F_{q^2}$, so that for all vectors $u,v$,
\[
(uM)^\phi U \left((vM)^{\phi q}\right)^\top=\lambda\cdot  uU(v^q)^\top
\]
where $\phi$ and the `power $q$' map act component-wise. If $u$ and $v$ are the canonical basis vectors $e_i$
and $e_j$, then this simplifies to
\[
e_i\left(M UM^{q\top}\right)e_j^\top=\left(\lambda \cdot e_iU(e_j^q)^\top\right)^{\phi^{-1}}=  e_i(\lambda^{\phi^{-1}} \cdot U)e_j^\top,
\]
or in other words, $M UM^{q\top}=\lambda^{\phi^{-1}}\cdot U$. 
Comparing the zero entries of $U$ with the entries of $M UM^{q\top}$ gives $d=i=a=e=0$.
Moreover, $bf^q+ch^q=0$. From Equation \eqref{eq:first}, $f=1-b$, $c=1-h$, and $\mu=1$.
Now the $(4,1)$-entry of $M UM^{q\top}$ is $\mu=1$, and so it follows that $\lambda=1$.
So we also have (from the diagonal entries), $bb^q+cc^q=1$, and $ff^q+hh^q=1$.
Therefore, $b+c=1$ (because $b(1-b^q)+c(1-c^q)=0$ and $bb^q+cc^q=1$)
and similarly, $f+h=1$. Hence, $M$ is of the form
\[
\begin{bmatrix}
    1 & 0 & 0 & 0\\
    0 & b & 1-b & 0\\
    0 & 1-b & b & 0\\
    0 & 0 & 0 & 1
\end{bmatrix}
\]
where $b+b^q=2bb^q$. 

The involutions
\[M_1:=\begin{bmatrix}
    0&0&0&1\\
    0&1&0&0\\
    0&0&1&0\\
    1&0&0&0
\end{bmatrix},\quad 
M_2:=\begin{bmatrix}
    0&0&0&1\\
    0&0&-1&-1\\
    0&-1&0&-1\\
    1&1&1&1
\end{bmatrix}
\]
centralise $C$, stabilise $\{P,Q,S\}$, and permute $\{P,Q,S\}$ (e.g., $M_1$ interchanges
$P$ and $Q$; $M_2$ acts as the cycle $(P \, S\, Q)$).
The group $K:=\langle M_1,M_2\rangle$ is isomorphic to $S_3$
and it follows that $D:=\langle K,C\rangle$ is the setwise stabiliser of $\{P,Q,S\}$ in $\PGammaU(4,q)$.
We have that $|S_3|=6$ and $\mathrm{Aut}(\F_{q^2})=2\log_p(q)$, and finally, we see that there are 
$q+1$ elements of $x\in \F_{q^2}$ such that $x+x^q=2x^{q+1}$: the splitting field of the polynomial $g(x)=2x^{q+1}-x^q-x=xh(x)$ is $\F_{q^2}$ and there are no multiple roots since $h'(x)=(q-1)x^{q-2}$ has only $0$ as a root. 
\end{proof}

\begin{lemma}\label{lemma:transitive}
The collineation group $G:=\PGammaU(4,q)$ of the Hermitian surface $\h(3,q^2)$ acts transitively on pairs of noncollinear points and transitively on triples of noncollinear points in perspective.
\end{lemma}

\begin{proof}
By \cite[Corollary 11.12]{Grove}, $\mathrm{PGU}(4,q)$ is transitive on the set of pairs of noncollinear points of
$\h(3,q^2)$. Consider the points $P(0,0,0,1)$, $Q(1,0,0,0)$, and $S(1,1,1,1)$. It is not difficult to see
that $[P,Q,S]=1$ and hence $\{P,Q,S\}$ are in perspective. 
By \cite[Lemma 4.2]{Bamberg:2021aa}, the number of triples in perspective is
$\frac{1}{6}q^5(q^2+1)(q^3+1)\cdot (q^3-q)(q-1)=\frac{1}{6}q^6(q^3+1)(q^4-1)(q-1)$.

By Lemma \ref{lemma:stabPQS}, the
setwise stabiliser in $\PGammaU(4,q)$ of $P,Q,S$ 
has size $12(q+1)f$ where $f=\log_p(q)$ and $p$ is the characteristic of $\F_q$.
Now $|\PGammaU(4,q)|=2q^6(q^3+1)(q^4-1)(q^2-1)f$ (see \cite[Table 2.1D]{kleidmanliebeck}) and so by the Orbit-Stabiliser Theorem, the size of the orbit of $\{P,Q,S\}$ under $\PGammaU(4,q)$ is $\frac{2q^6(q^3+1)(q^4-1)(q^2-1)f}{12(q+1)f}=\frac{1}{6}q^6(q^3+1)(q^4-1)(q-1)$; the number
of perspective triples.
\end{proof}

As a corollary, the last two results (Lemma \ref{lemma:stabPQS} and Lemma \ref{lemma:transitive}) give the following.

\begin{lemma}\label{classical}
Let $\s$ be a set of points of $\h(3,q^2)$ containing $P,Q,S$ such that there is an element of the collineation
group $\PGammaU(4,q)$ of $\h(3,q^2)$ mapping $\mathcal{V}$ onto $\s$. Then $\s$ is of the form
\[
\left\{ \left(1,\; (1-x)a^q+xa,\; x a^q+(1-x)a,\; a^{q+1}\right):a\in\F_{q^2}\}\cup\{(0,0,0,1)\right\}
\]
for some fixed $x\in \F_{q^2}$ satisfying $x+x^q=2x^{q+1}$, and conversely.
\end{lemma}



\begin{lemma} \label{lemma:elliptic}
Let $\mathcal{E}=\{(1,a,b,\frac{a^2+b^2}{2})\mid a,b\in \F_q \}$ be an elliptic quadric contained in a subgeometry $\PG(3,q)$ of $\PG(3,q^2)$. Suppose that $q\equiv 3\pmod{4}$, then $\mathcal{E}$ is equivalent under the collineation group of $\h(3,q^2)$ equivalent to the Hermitian Veronesean.  
\end{lemma}

\begin{proof}
We first show that there is an element $x\in \F_{q^2}\setminus \F_q$ with $2\norm(x)=\trace(x)=1.$
Recall that an equation of the form $X^2-tX+n=0$, where $t,n\in \F_q$ has either roots in $\F_q$, or has roots $\gamma,\gamma^q\in \F_{q^2}$ with $\trace(\gamma)=t$ and $\norm(\gamma)=n$, and vice versa, $\gamma,\gamma^q\in \F_{q^2}$ are the roots of $X^2-\trace(\gamma)X+\norm(\gamma)=0$. Using that $\trace(x)=1$ and $\norm(x)=\frac{1}{2}$, we see that $x$ would need to be the root of the equation $X^2-X+\frac{1}{2}=0.$
This equation has discriminant $-1$, which is a non-square if $q\equiv 3\pmod{4}$, so we find a solution $x\in \F_{q^2}\setminus \F_q$ with $2\norm(x)=1=\trace(x)$.

Fix this element $x$. We will now show that any point of the form $(1,a,b,\frac{a^2+b^2}{2})$, $a,b\in \F_q$, can be written as $(1,(1-x)\alpha^q+x\alpha,x\alpha^q+(1-x)\alpha,\alpha^{q+1})$ for some $\alpha\in \F_{q^2}$.
Note that, since $x^q=1-x$, we have that $(1-x)\alpha^q+x\alpha=\trace(x\alpha)$ and $x\alpha^q+(1-x)\alpha=\trace(x^q\alpha)$.
Consider the map $\psi:\alpha\to (\trace(\alpha x),\trace(\alpha^q x))$, which is an $\F_q$-linear map from $\F_{q^2}$ to $\F_q^2$. We will show that $\psi$ is injective (and therefore, surjective) by showing that $\ker(\psi)=\{0\}$. So let $\alpha\in \F_{q^2}$ with $\trace(\alpha x)=\alpha x+\alpha^q x^q=0$ and $\trace(\alpha x^q)=\alpha x^q+\alpha ^q x=0$, then it follows that $(\alpha-\alpha^q)(x-x^q)=0$. Since $x\notin \F_q$, $x-x^q=0$, and it follows that $\alpha\in \F_q$. When $\alpha\in \F_q$, we have that $0=\trace(\alpha x)=\alpha \trace(x)=\alpha$. Hence, our claim follows.
\end{proof}

We will now prove our first main result -- Theorem \ref{main1}.

\begin{theorem} 
Let $\s$ be a set of mutually noncollinear points on $\h(3,q^2)$, $q$ odd.
Suppose there are four non-coplanar points $P,Q_1,Q_2,Q_3\in \s$, such that all triangles $PQ_iR$, $i=1,2,3$, $R\in \s\setminus \{P,Q_i\}$, are in perspective.
Then $|\s|\le q^2+1$ and equality holds if and only if $\s$ is equivalent under the collineation group of $\h(3,q^2)$ to the Hermitian Veronesean $\mathcal{V}$.
\end{theorem}

\begin{proof}
    Since the stabiliser of $\h(3,q^2)$ is transitive on triples in perspective by Lemma \ref{lemma:transitive}, we can pick $P(0,0,0,1)$, $Q_1(1,0,0,0)$, and $Q_2(1,1,1,1)$.
    Now let $R(1,a,b,c)$ be a point of $\s$, different from $P,Q_1,Q_2$; the first coordinate of $R$ can be taken as $1$ since $R$ is not collinear with $P$. Since $R\in \h(3,q^2)$, we find that $\trace(c)=\norm(a)+\norm(b)$. Since $PQ_1R$ is in perspective, the Segre invariant $[P,Q_1,R]=c^q\in \F_q$ which implies that $c\in \F_q$. Since $PQ_2R$ is in perspective, the Segre invariant $[P,Q_2,R]=1-a^q-b^q+c^q\in \F_q$ which forces $a+b\in \F_q$. It follows that $2c=\trace(c)=N(a)+N(b).$

Let the fourth point $Q_3$ have coordinates $Q_3(1,a,b,c)$ and consider a basis $\{1,\omega\}$ with $\trace(\omega)=0$. Since $a+b\in \F_q$, we find that  $a=a_0+a_1\omega$ and $b=b_0-a_1\omega$.

Consider now an arbitrary element $\tilde{R}(1,\tilde{a},\tilde{b},\tilde{c})$ in $\s$. We have seen that $\tilde{a}+\tilde{b}\in \F_q$ and $\tilde{c}\in \F_q$.
Write $\tilde{a}=\tilde{a_0}+\tilde{a_1}\omega$ and $\tilde{b}=\tilde{b_0}-\tilde{a_1}\omega$.
Since $PQ_3\tilde{R}$ is in perspective, we find that the Segre invariant $[P,Q_3,\tilde{R}]=c+\tilde{c}^q-a\tilde{a}^q-b\tilde{b}^q\in \F_q
$. 
It follows that $(a_0+a_1\omega)(\tilde{a_0}+\tilde{a_1}\omega^q)+(b_0-a_1\omega)(\tilde{b_0}-\tilde{a_1}\omega^q)\in \F_q$.

Using that $\omega^q=-\omega$ and putting the coefficient of $\omega$ in this expression equal to zero yields 
\begin{align}\tilde{a_1}(b_0-a_0)=a_1(\tilde{b_0}-\tilde{a_0}).\label{eq:a1}\end{align}
First assume that $a\in \F_q$, that is, that $a_1=0$. 
Note that in this case, necessarily $a_0\neq b_0$, since otherwise the points $P,Q_1,Q_2,Q_3$ would all be contained in the plane $Y=Z$, a contradiction. Equation \eqref{eq:a1} then implies that $\tilde{a_1}=0$, that is, $\tilde{a}\in \F_q$. It follows that every point $R$ in $\s$, different from $P$, has coordinates $(1,\tilde{a},\tilde{b},\frac{\tilde{a}^2+\tilde{b^2}}{2}),$ where $\tilde{a},\tilde{b}$ are in $\F_q$. It then immediately follows that $|\s|\le q^2+1$ and that equality holds if and only if $\s$ is the point set of an elliptic quadric contained in a subgeometry $\PG(3,q)$. If $q\equiv 3\pmod 4$, then we know from Lemma \ref{lemma:elliptic} that $\s$ is equivalent to a Hermitian Veronesean. 
If $q\equiv 1 \pmod 4$, then $-1$ is a square, so there exist $a,b \in \F_q$ such that $(a,b)\ne (0,0)$ and $a^2+b^2=0$. It follows that $R(1,a,b,0)$ and $Q_1(1,0,0,0)$ are collinear on $\h(3,q^2)$, a contradiction. So the elliptic quadric contained in as subgeometry does not form an example of set $\s$ satisfying our conditions when $q\equiv 1 \pmod{4}.$

We may now assume that $a_1\neq 0$. Put $\lambda=\frac{b_0-a_0}{a_1}$. We see that \eqref{eq:a1} implies that $\tilde{b_0}=\lambda \tilde{a_1}+\tilde{a_0}$, and hence, that all points $\tilde{R}\in \s$, different from $P$, are of the form 
$\tilde{R}(1,\tilde{a_0}+\tilde{a_1}\omega, \lambda \tilde{a_1}+\tilde{a_0}-\tilde{a_1}\omega,\tilde{c}),$ where $2\tilde{c}=\norm(\tilde{a_0}+\tilde{a_1}\omega)+\norm(\lambda \tilde{a_1}+\tilde{a_0}-\tilde{a_1}\omega)$.
In particular, we see that there are at most $q^2$ values for $a_0,a_1$, and hence, at most $q^2$ points in $\s$, different from $(0,0,0,1)$.
Suppose now that equality holds, so $|\s|=q^2+1$, then we have that every two elements $\tilde{a_0}, \tilde{a_1}$ determine a point $\tilde{R}(1,\tilde{a_0}+\tilde{a_1}\omega, \lambda \tilde{a_1}+\tilde{a_0}-\tilde{a_1}\omega,\tilde{c}),$ where $2\tilde{c}=\norm(\tilde{a_0}+\tilde{a_1}\omega)+\norm(\lambda \tilde{a_1}+\tilde{a_0}-\tilde{a_1}\omega)$ of $\s.$
The set $\s$ consists of mutually noncollinear points, so in particular, no point in $ \s$ should be collinear with $Q_1(1,0,0,0)$ (apart from $Q_1$, of course). A point $\tilde{R}$ is collinear with $Q_1$ if and only if $\tilde{c}=0$, which happens if and only if $\norm(\tilde{a_0}+\tilde{a_1}\omega)+\norm(\tilde{a_0}+\lambda \tilde{a_1}-\tilde{a_1}\omega)=0$.  This is $(\lambda^2-2\omega^2)\tilde{a_1}^2+2\lambda \tilde{a_0}\tilde{a_1}+2\tilde{a_0}^2=0,$
which is a quadratic equation
    in $\tilde{a_1}$ with coefficients in $\F_q$ (note that $\omega^2=-\omega^{q+1}=-\norm(\omega)$). This equation has discriminant
    \[    \Delta=4\tilde{a_0}^2 \left(4\omega^2- \lambda^2\right).
    \]
    Since this discriminant should be a non-square for every value of $\tilde{a_0}$, we see that $4\omega^2-\lambda^2$ is a non-square in $\F_q$.

We will now show that
\[
(\tilde{a_0}+\tilde{a_1}\omega, \lambda \tilde{a_1}+\tilde{a_0}-\tilde{a_1}\omega)=
\left((1-x)\alpha^q+x\alpha,x \alpha^q+(1-x)\alpha\right)
\]
for some $\alpha,x\in\F_{q^2}$ with $x+x^q=2x^{q+1}$. Lemma \ref{classical} then shows that $\s$ is 
equivalent under the collineation group of $\h(3,q^2)$ to the Hermitian Veronesean.

Since $4\omega^2-\lambda^2$ is a non-square in $\F_q$, it is a square in $\F_{q^2}$, say $4\omega^2-\lambda^2=\sigma^2$, where $\sigma\notin \F_q$. Now take
\[\alpha:=\tilde{a_0} +\tfrac12 \tilde{a_1} (\lambda+\sigma^q),\quad 
x:=\frac{\lambda -  2\omega}{\sigma-\sigma^q}+\frac{1}{2}.
\]
Since $\sigma\notin \F_q$ but $\sigma^2\in \F_q$, we know that $\trace(\sigma)=0$. It is now a straightforward calculation to check that indeed, $\alpha^q+x(\alpha-\alpha^q)=\tilde{a_0}+\tilde{a_1}\omega$, $\alpha+x(\alpha^q-\alpha)=\tilde{a_0}+\lambda \tilde{a_1}-\tilde{a_1}\omega$ and $2\norm(x)=\trace(x)$
\end{proof}

\begin{remark}
    The condition that the points $P,Q_1,Q_2,Q_3$ are not coplanar is necessary. If $Q_3$ is contained in the plane spanned by $P(0,0,0,1)$, $Q_1(1,0,0,0)$, and $Q_2(1,1,1,1)$, then it is not too hard to see that $Q_3(1,a,a,a^2)$ for some $a\in \F_q$. It easily follows from the arguments of our proof that the set $\s'$ of points of the form $\tilde{R}(1,\tilde{a},\lambda-\tilde{a},\frac{2\norm(\tilde{a})-\lambda\trace(\tilde{a})+\lambda^2}{2})$, where $\tilde{a}\in \F_{q^2}$ and $\lambda\in \F_q$, is a set of $q^2\cdot q$ points such that the Segre invariants $[P,Q_i,\tilde{R}]\in \F_q$ for $i=1,2,3$. For each choice of $\tilde{a}$, at most $2$ choices of $\lambda$ can give rise to a triangle with $[P,Q_1,\tilde{R}]=0$, since expressing that $[P,Q_1,\tilde{R}]=\tilde{c}=0$ gives a quadratic equation in $\lambda$. Similarly, for each value of $\tilde{a}$, we find at most $2$ values for $\lambda$ such that $\tilde{R}$ spans a degenerate plane with $P,Q_2$ or with $P,Q_3$. So we find at least $q^2(q-6)$ points different from $P$ in $\s'$, which is larger than $q^2$ if $q\ge 7$.
 \end{remark}

\section{A characterisation via sublines}

Let $\h(3,q^2)$, $q$ odd, be defined by the Hermitian form $h(X,Y):=X_1Y_4^q+X_4Y_1^q-X_2Y_2^q-X_3Y_3^q$. Special sets of $\h(3,q^2)$, in general, have an interesting property:
every point $Z$ not in $\s$ is collinear with $0$ or $2$ elements of $\s$. In particular, for every point $P$ in $\s$ and every totally isotropic line $\ell$ of $\h(3,q^2)$ through $P$, it holds that every point of $\ell$ is collinear (in $\h(3,q^2)$) with exactly one point of $\s$, different from $P$. We will replace the condition that $\s$ is a special set with a local condition where we fix a point $P$ and a totally isotropic line $\ell$ through $P$.

\begin{lemma}\label{lemma:bijection}
    Suppose that a set $\s$ of points on $\h(3,q^2)$ contains the point $P(0,0,0,1)$, and let
    $\ell$ be the span of $P$ with $(0,1,\omega,0)$, where $\omega^{q+1}+1=0$. (So $\ell$ is totally isotropic). Suppose that every point on $\ell$ is collinear with exactly one point of $\s\setminus\{P\}$. Then there is a natural bijection $F_\ell$ between the points of $\ell$
    and the points of $\s$ defined by the following: $F_\ell(P)=P$,
    and for any point $Y$ on $\ell$ not equal to $P$, we let $F_\ell(Y)$ be the unique
    point of $\s\backslash\{P\}$ collinear with $Y$.
\end{lemma}

\begin{corollary}
 Let $\s$ be as in Lemma \ref{lemma:bijection}.
There exists a function $f:\F_{q^2}\to \F_{q^2}$ such that
every point of $\s\backslash\{P\}$ is of the form $(1,f(t)+t,(f(t)-t)\omega,c_t)$ where $t\in\F_{q^2}$,
and $c_t$ is some element of $\F_{q^2}$ such that $\trace(c_t)=2 \trace(f(t)t^q)$.
\end{corollary}

\begin{proof}
Let $\ell$ be the span of $P$ with $(0,1,\omega,0)$, where $\omega^{q+1}+1=0$. (So $\ell$ is totally isotropic). Every point $Y$ on $\ell$, not equal to $P$, is of the form \[(0,1,\omega,2t^q),\quad t\in \mathbb{F}_{q^2}\] and $Y$ is collinear (in $\h(3,q^2)$) to a unique element of $\s\backslash\{P\}$,
which we denote by $F_\ell(Y)$ in Lemma \ref{lemma:bijection}. 
Let $F_\ell(Y):=(1,a_t,b_t,c_t)$. Write $a_t=t+f(t)$ (this is always possible), then $(0,1,\omega,2t^q)$ and $(1,t+f(t),b_t,c_t)$ being collinear (in $\h(3,q^2)$) implies that 
$t-b_t\omega^q -f(t)=0$, and hence, that $b_t=(f(t)-t)\omega$. So we
can write this point as $(1,f(t)+t,(f(t)-t)\omega,c_t)$. 
Moreover, since this point is totally isotropic, we have 
\[
0=\trace(c_t)-(f(t)+t)^{q+1}-(f(t)-t)^{q+1}\omega^{q+1}=\trace(c_t)-2 \trace(f(t)t^q).\qedhere
\]
\end{proof}

\begin{corollary}\label{fzero}
     Let $\s$ be as in Lemma \ref{lemma:bijection} and suppose that $\s$ contains $Q(1,0,0,0)$. Then $f(0)=0$ and $c_0=0$.
\end{corollary}

\begin{proof}
    Since $\h(3,q^2)$ is a generalised quadrangle, the point $Q\in \s$ is collinear in $\h(3,q^2)$ with a unique element on $\ell$, which is different from $P$. Hence, $Q=F_\ell(Y)$ for some $Y\in \ell$ with $Y(0,1,\omega,2t^q)$. Since $t+f(t)=0$ and $(f(t)-t)\omega=0$, it follows that $t=f(t)=0$, and so $f(0)=0$. It now follows that $c_0=0$.
\end{proof}



\begin{lemma}\label{lemma:flinear} 
Let $\s$ be as in Lemma \ref{lemma:bijection} and assume that $Q(1,0,0,0)\in \s.$
Suppose that for every Baer subline $b$ of $\ell$ containing $P$,
the image of $b$ under $F_\ell$ lies in a plane. Then $f$ is $\F_q$-linear.
Conversely, if $f$ is $\F_q$-linear then the set $\s=\{(1,f(t)+t,(f(t)-t)\omega,\trace(f(t)t^q))\mid t\in \F_{q^2}\}$ is a set of $q^2$ points on $\h(3,q^2)$ such that for every Baer subline $b$ of $\ell$ containing $P$, all points of $\s$ collinear with a point of $b$ are coplanar.
\end{lemma}

\begin{proof} 
First assume that every Baer subline $b$ of $\ell$ containing $P$, the image of $b$ under $F_\ell$ lies in a plane. Let $u_1,u_2$ be different elements of $\F_{q^2}$. Then
the Baer subline\footnote{Write $B=A+\alpha P$. Then the Baer subline on $\{A,P,B\} $ consists of points of the form $A+\lambda\alpha\cdot P$ and the point $P$. Let $A=(0,1,\omega,u_1)$ and $B=(0,1,\omega,u_2)$. Now $\alpha=u_2-u_1$ and so $A+\lambda\alpha\cdot P=(0,1,\omega,u_1+\lambda(u_2-u_1))$.} $b$ on $P$, $(0,1,\omega,u_1)$, and $(0,1,\omega,u_2)$
consists of $P$ together with the points of the form
\[
(0,1,\omega, u_1+\lambda(u_2-u_1))
\]
where $\lambda\in\F_q$. Consider $P$ and the image of $(0,1,\omega,u_1)$, $(0,1,\omega,u_2)$,
and $(0,1,\omega, u_1+\lambda(u_2-u_1))$ under $F_\ell$, written
as the rows of a matrix $A$:
\[
A=\begin{bmatrix}
&0,&0,&0,&1\\
& 1, &f(u_1)+u_1 ,  &(f(u_1)-u_1)\omega, & * \\
& 1 , &f(u_2)+u_2 , &(f(u_2)-u_2)\omega, & * \\
& 1 , &f( \lambda (u_2-u_1)+u_1)+\lambda (u_2-u_1)+u_1 , 
&(f(\lambda (u_2-u_1)+u_1)-\lambda (u_2-u_1)-u_1)\omega, &* \\
\end{bmatrix},
\]
where we have suppressed the last column. The determinant of this
matrix is
\(
2 \omega (u_1-u_2) (f(-\lambda u_1+\lambda u_2+u_1)+(\lambda-1) f(u_1)-\lambda f(u_2))
\)
which is zero if and only if
\[
f((1-\lambda) u_1+\lambda u_2)=(1-\lambda) f(u_1)+\lambda f(u_2).
\]
Our assumption implies that this holds for all $u_1\ne u_2$ and all $\lambda\in\F_q$.
Taking $u_1=0$ and using Corollary \ref{fzero}, implies that \(f(\lambda u_2)=\lambda f(u_2)\)
for all $\lambda\in\F_q$, and all $u_2\in \F_{q^2}^*$. Now
taking $\lambda=\tfrac12$ yields $f(\tfrac12( u_1+u_2))=\tfrac{1}{2}(f(u_1)+f(u_2))$
and so by our first property, \(f(u_1+u_2)=f(u_1)+f(u_2)\)
for all $u_1\neq u_2\in \F_{q^2}$. Finally, if $u_1=u_2$, then $f(u_1+u_2)=f(2u_1)=2f(u_1)=f(u_1)+f(u_1)=f(u_1)+f(u_2),$ so we conclude that $f(\lambda u)=\lambda f(u)$ for all $u\in \F_{q^2}$ and $\lambda \in \F_q$ and $f(u_1+u_2)=f(u_1)+f(u_2)$ for all $u_1,u_2\in \F_{q^2}$, that is, $f$ is $\F_q$-linear.

Conversely, assume that $f$ is $\F_q$-linear. A point on $\ell$, different from $P$, is given by $(0,1,\omega,2t^q)$ and we see that this point is collinear with the point $(1,f(t)+t,(f(t)-t)\omega,\trace{f(t)t^q})$. Moreover, $(1,f(t)+t,(f(t)-t)\omega,\trace{f(t)t^q})$ is the unique such point in $\s$, and it is contained in $\h(3,q^2)$. As deduced above, every subline on $\ell$ through the point $P$, consists of points of the form $(0,1,\omega,u_1)$, $(0,1,\omega,u_2)$, $(0,1,\omega,u_1+\lambda(u_2-u_1))$. The points of $\s$ collinear with those points are given by the rows of the matrix $A$, whose determinant is $0$ because when $f$ is $\F_q$-linear, $f((1-\lambda)u_1+\lambda u_2)=(1-\lambda)f(u_1)+\lambda f(u_2)$.
\end{proof}

Under the assumptions of Lemma \ref{lemma:flinear}, we know that $f$ is an $\F_q$-linear map on $\F_{q^2}.$ So we can write $f$ as a linearised polynomial: $f(t)=\alpha t+\beta t^q$ where $\alpha,\beta$ are fixed elements in  $ \F_{q^2}$.
We then see that the points of $\s\backslash\{P\}$, consists of elements of the form
\[
(1,t+\alpha t+\beta t^q,(\alpha t+\beta t^q-t)\omega, c_t),
\]
for some $t\in \F_{q^2}.$

The map $f$ featured in the results above, enjoys some interesting properties.

\begin{proposition}
Let $\s$ be as in Lemma \ref{lemma:bijection}.
\begin{enumerate}[(a)]
    \item  Suppose for every Baer subline $b$ of $\ell$ containing $P$,
    the image of $b$ under $F_\ell$ lies in a plane. Then
    $\langle P, R_1, R_2\rangle$ is a nondegenerate plane if and only if
    $\trace(f(t_1-t_2)(t_1-t_2)^q)\ne 0$.
    \item Suppose $PQR$ is in perspective for all $R\in \s\backslash\{P,Q\}$. Then
    $\langle P,R_1,R_2\rangle$ is a nondegenerate plane, for all distinct $R_1,R_2\in \s\backslash\{P\}$.
\end{enumerate}
\end{proposition}

\begin{proof}\leavevmode
\begin{enumerate}[(a)]
    \item 
    We know $\trace(c_t)=2 \trace(f(t)t^q)$. 
    Also, since $f$ is $\F_q$-linear,
\begin{align*}    
    \trace[P,R_1,R_2]&=\trace(c_{t_1}+c_{t_2})-2\trace(f(t_2)t_1^q+f(t_1)t_2^q)\\
    &=2 \trace(f(t_1)t_1^q+f(t_2)t_2^q)-2\trace(f(t_2)t_1^q+f(t_1)t_2^q)\\
    &=2 \trace(f(t_1-t_2)(t_1-t_2)^q)).
\end{align*}   
    So $\langle P, R_1, R_2\rangle$ is nondegenerate if and only if $\trace(f(t_1-t_2)(t_1-t_2)^q)\ne 0$.
    \item Recall that $c_{t_1-t_2}\ne 0$, since $t_1-t_2\ne 0$.
Now
\([P,R_1,R_2]=-2 t_1 f(t_2)^q-2 t_2^q f(t_1)+\trace(t_1 f(t_1)^q+t_2 f(t_2)^q)\)
and so
\begin{align*}
\trace([P,R_1,R_2])&=2 \cdot \trace\left( \left(f(t_1)-f(t_2)\right)(t_1-t_2)^q\right)\\
&=2 \left(\trace(\alpha) (t_{1}-t_{2})^{q+1}+\trace(\beta^q (t_{1}-t_{2})^2)\right)\\
&=2c_{t_1-t_2}\ne 0.\qedhere
\end{align*}
\end{enumerate}
\end{proof}



\begin{lemma}\label{lemma:ct} 
Let $P(0,0,0,1)$, $Q(1,0,0,0)$ and $R(1,t+\alpha t+\beta t^q,(\alpha t+\beta t^q-t)\omega, c_t)$ be points of $\h(3,q^2)$. Then $PQR$ is in perspective if and only if $c_t=t^{q+1}\trace(\alpha)+\trace(t^2\beta)$ and $c_t\neq 0.$
\end{lemma}
\begin{proof}
In order for $R$  to be totally isotropic, we must have
\[
\trace(c_t)=2  \trace(t^{q+1}\alpha+t^2 \beta^q).
\]
Now $[P,Q,R]=c_t^q$ and hence $\trace(c_t)\ne 0$ because $P$, $Q$, and $R$ span a nondegenerate plane. So $[P,Q,R]\in \F_q^*$ and hence $c_t\in \F_q^*$.
This means
$c_t=t^{q+1}\trace(\alpha)+\trace(t^2 \beta^q)$.
\end{proof}

\begin{defn}\label{defn:salphabeta}
    Let \[\s_{\alpha,\beta}=\{P(0,0,0,1)\}\cup\{(1,t+\alpha t+\beta t^q,(\alpha t+\beta t^q-t)\omega,t^{q+1}\trace(\alpha)+\trace(t^2\beta))\mid t\in \F_{q^2}\}.\]
    Note that $Q(1,0,0,0)\in \s_{\alpha,\beta}.$ Furthermore, let 
    \[
    R_i:=(1,(\alpha+1) t_i+\beta t_i^q,((\alpha-1)t_i+\beta t_i^q)\omega, t_i^{q+1}\trace(\alpha)+\trace(t_i^2 \beta^q)).
    \]
\end{defn}

\begin{corollary} 
    Suppose $PQR$ is in perspective for all $R\in \s_{\alpha,\beta}\backslash\{P,Q\}$. Then $\trace(\alpha+\beta)\neq 0.$ 
\end{corollary}
\begin{proof}
    We have seen in Lemma \ref{lemma:ct} that $c_t\neq 0$ for all $t$, in particular, for $t=1$, we see that $\trace(\alpha+\beta)=\trace(\alpha)+\trace(\beta)\neq 0.$
\end{proof}



The following is a straight-forward calculation:

\begin{lemma}\label{PR1R2}
    Consider $R_1,R_2\in \s_{\alpha,\beta}$ with $t_1\neq t_2$. Then 
    \[
    [P,R_1,R_2]=\trace(\alpha) \left( t_1^{q+1}-2 t_1t_2^q+t_2^{q+1}\right)+\trace( (t_1-t_2)^2\beta^q).
    \]
    
    Moreover, $[P,R_1,R_2]\in \F_q$ if and only if $\trace(\alpha)=0$ or $t_1t_2^q\in \F_q$.
\end{lemma}

\begin{lemma}\label{lemma:classical}
 Suppose $PQR$ is in perspective for all $R\in \s_{\alpha,\beta}\backslash\{P,Q\}$. Then
    $\s_{\alpha,\beta}$ is equivalent under the collineation group of $\h(3,q^2)$ to the Hermitian Veronesean $\mathcal{V}$ if and only if $\trace(\alpha)=0$.
\end{lemma}

\begin{proof}
Suppose $\s$ is equivalent under the collineation group of $\h(3,q^2)$ to the Hermitian Veronesean $\mathcal{V}$. Then $[P,R_1,R_2]\in \F_q$ for all $R_1,R_2\in \s\backslash\{P\}$.
By Lemma \ref{PR1R2}, we have $\trace(\alpha)=0$, or $t_1t_2^q\in \F_q$ for every distinct
pair $t_1,t_2\in \F_{q^2}$. Clearly the latter does not hold (e.g., take $t_2=1$ and $t_1\notin\F_q$),
and so $\trace(\alpha)=0$.

Conversely, suppose $\trace(\alpha)=0$.
Then, by Lemma \ref{PR1R2}, $[P,R_1,R_2]=\trace( (t_1-t_2)^2\beta^q)$ for every distinct pair $R_1,R_2\in\s_{\alpha,\beta}\backslash\{P\}$. So each $[P,R_1,R_2]$ lies in $\F_q^\ast$, meaning that all triangles on $P$ are in perspective. Corollary \ref{justonepoint} shows that $\s_{\alpha,\beta}$ is  equivalent under the collineation group of $\h(3,q^2)$ to the Hermitian Veronesean $\mathcal{V}$. 
\end{proof}

\begin{corollary}
    Suppose some $PR_1R_2$ is in perspective, where the $R_1,R_2\in \s_{\alpha,\beta}$
    correspond to $t_1,t_2$ such that
    $t_2t_1^q\ne t_1t_2^q$. Then $\s_{\alpha,\beta}$ is equivalent under the collineation group of $\h(3,q^2)$ to the Hermitian Veronesean $\mathcal{V}$.
\end{corollary}

\begin{proof}
Notice that 
\[
[P,R_1,R_2]=(t_2t_1^q-t_1t_2^q)\trace(\alpha)+\norm(t_1-t_2)\trace(\alpha)+\trace((t_1-t_2)^2\beta^q)
\]
and so $(t_2t_1^q-t_1t_2^q)\trace(\alpha)\in\F_q$.
Now we assumed $t_2t_1^q\ne t_1t_2^q$, which is equivalent to $t_2t_1^q-t_1t_2^q\notin \F_q$,
and hence $\trace(\alpha)=0$. 
\end{proof}

\begin{remark} The set
$\s_{\alpha,\beta}$ is a set of $q^2$ points such that $[P,Q,R]\in \F_q$ for all points $R\in \s$ and such that for every subline $b$ in $\ell$, the set of points in $\s$ collinear with the points of $b$ is coplanar. It is possible to find $\alpha_0,\beta_0$ with $\trace(\alpha_0)\neq 0$ such that $[P,Q,R]\in \F_q^*$ for all $R\in \s_{\alpha_0,\beta_0}\setminus\{P,Q\}$, and hence $\s_{\alpha_0,\beta_0}$ is a set of points such that $PQR$ is in perspective for all $R$ but it is not a special set.
To find such $\alpha_0,\beta_0$ we need to express that the equation $[P,Q,R]=\norm(t)\trace(\alpha)+\trace(\beta^q t^2)=0$ only has the solution $t=0$. Writing $\alpha_0=c_0+c_1\omega$, $\beta_0=d_0+d_1\omega$ and $t=t_0+t_1\omega$ for some $c_i,d_i,t_i\in \F_q$ and $\omega\in\F_{q^2}$ with $\trace(\omega)=0,$ we find that any $c_0\neq 0,c_1,d_0,d_1$ such that $d_1^2\norm(\omega)^2-\norm(\omega)(d_0^2-c_0^2)$ is not a square in $\F_q$ yields an admissible such $\alpha_0,\beta_0$. 
\end{remark}

\begin{lemma} Suppose $t_2,t_3\in \F_q$, then 
\begin{align*}
&\trace([R_1,R_2,R_3])=(t_2-t_3)^2 \trace(\alpha+\beta)\cdot\\ &\qquad\trace\left(\left( (t_1^{q+1}-2 t_2t_1^q+t_2^2) \trace(\alpha)+\trace(\beta^q (t_1-t_2)^2)\right) \left( 
\left( t_1^{q+1}-2 t_1t_3+t_3^2\right)\trace(\alpha)+\trace(\beta^q (t_1-t_3)^2)\right)\right).
\end{align*}
\end{lemma}
\begin{proof}
    Note that $[R_1,R_2,R_3]=[P,R_1,R_2][P,R_2,R_3][P,R_3,R_1].$ 
    A direct (yet tedious) calculation now yields the result.
\end{proof}

\begin{lemma}\label{lemma:tracenonzero}
    Suppose $\trace(\alpha)\ne 0$. 
    Then there exist $t_1\in\F_{q^2}$, $t_2,t_3\in \F_q$, such that
    $\trace[R_1,R_2,R_3]=0$ and $(t_1,t_2,t_3)$ is a distinct triple.
\end{lemma}

\begin{proof}
    Consider a basis $\{1,\xi\}$ for $\F_{q^2}$ with $\trace(\xi)=0$. Then $\xi^q=-\xi$ and we denote $\norm(\xi)=-\xi^2$ by $\kappa$. Note that $\kappa\in \F_q^\ast$ and that $-\kappa=\xi^2$ is a non-square in $\F_q$, since $\xi\in \F_{q^2}\setminus \F_q$. Let $t_1=y_0+y_1\xi$, then $\norm(t_1)=y_0^2+\kappa y_1^2$.
    Let $\beta^q=b_0+b_1\xi$ and $\trace(\alpha)=a$. 
    We see that, for $t_2,t_3\in \F_q$, $t_1^{q+1}-2t_2t_1^q+t_2^2=(y_0-t_2)^2+\kappa y_1^2+2y_1t_2\xi$ and that $\trace((b_0+b_1\xi)(t_1-t_2)^2)=2b_0(y_0-t_2)^2-2b_0\kappa y_1^2-4b_1y_1\kappa(y_0-t_2)$. Using $s_1=y_0-t_2$, we can now rewrite $ (t_1^{q+1}-2 t_2t_1^q+t_2^2) \trace(\alpha)+\trace(\beta^q (t_1-t_2)^2)$ as
    \[a(s_1^2+\kappa y_1^2+2y_1t_2\xi)+2b_0(s_1^2-\kappa y_1^2) -4b_1y_1\kappa s_1.\] Similarly, with $s_2=y_0-t_3$ we rewrite $ (t_1^{q+1}-2 t_3t_1^q+t_3^2) \trace(\alpha)+\trace(\beta^q (t_1-t_3)^2)$ as
    $a(s_2^2+\kappa y_1^2+2y_1t_3\xi)+2b_0(s_2^2-\kappa y_1^2) -4b_1y_1\kappa s_2$.
If $s_1\neq s_2$ (that is, if $t_2\neq t_3$), then $\trace[R_1,R_2,R_3]=0$ if and only if \[
 \trace\left(\left( (t_1^{q+1}-2 t_2t_1^q+t_2^2) \trace(\alpha)+\trace(\beta^q (t_1-t_2)^2)\right) \left( 
\left( t_1^{q+1}-2 t_1t_3+t_3^2\right)\trace(\alpha)+\trace(\beta^q (t_1-t_3)^2)\right)\right)
=0.\]
Using the expressions deduced earlier, together with $\trace(\xi)=0, \xi^2=-\kappa$ and $t_3=s_1-s_2+t_2$, this in turn is equivalent to
\begin{align}(a(s_1^2+\kappa y_1^2)+2b_0(s_1^2-\kappa y_1^2) -4b_1y_1\kappa s_1)(a(s_2^2+\kappa y_1^2)+2b_0(s_2^2-\kappa y_1^2) -4b_1y_1\kappa s_2)=\nonumber\\4y_1^2t_2(s_1-s_2+t_2)\kappa\nonumber\\
\iff(s_1^2(a+2b_0)+\kappa y_1^2(a-2b_0) -4b_1y_1\kappa s_1)(s_2^2(a+2b_0)+\kappa y_1^2(a-2b_0) -4b_1y_1\kappa s_2)\nonumber=\\4y_1^2t_2(s_1-s_2+t_2)\kappa
\label{eq:s1}.\end{align}
We will show that there is a solution $(s_1,s_2,y_1,t_2)\in \F_q^4$ of \eqref{eq:s1}, where $y_1\neq 0$ and $s_1\neq s_2$ since any such solution will give us $t_1\in \F_q^2\setminus\{\F_q\}$ and $t_2\neq t_3\in \F_q$.
If one of the quadratic factors on the left hand side of the equation has nonzero roots, say $(s_1=\tau,y_1=\rho)\neq (0,0)$ then we can take $(\tau,s_2,\rho,0)$ where $s_2\neq \tau$ is arbitrary. So now assume that there are no $s,y\in \F_q^*$ such that $s^2(a+2b_0)+\kappa y^2(a-2b_0)-4b_1\kappa ys=0.$ Using that $4b_1^2\kappa^2-4b_0^2\kappa=\kappa \norm(\beta)$ we see that this assumption is equivalent to the assumption that $\delta:=\kappa(4\norm(\beta)-a^2)$ is a non-square.

{\underline{ Case 1: $q\equiv 3\pmod{4}$}}.  

Let $t_2=a$ and $s_2=0$, Then, after dividing by $\kappa y_1^2$, we find a quadratic equation in $y_1$ with discriminant
\[(a-2b_0)^2(s_1^2(4b_0^2\kappa+4b_1\kappa^2b_1^2-a^2\kappa)+4a\kappa(s_1+a)).\]

This is a square if there exists some $s_1\neq 0$ such that \[\kappa(s_1^2(4\norm(\beta)-a^2)+4as_1+4a^2)\] is a square. Note that $4\norm(\beta)-a^2\neq 0$ since we have assumed that $\delta$ is not a square, so we can consider $s_1=\frac{4a}{a^2-4 \norm(\beta)}$. Since $a\neq 0$, $s_1\neq 0$, and  $\kappa(s_1^2(4\norm(\beta)-a^2)+4as_1+4a^2)=4a^2\kappa$, which is a square since  $\kappa$ is a square: since $q\equiv 3\pmod{4}$, $-1$ is not a square and $-\kappa$ is not a square. 

{\underline{ Case 2: $q\equiv 1\pmod{4}$}}.


First assume that $b_1\neq 0$.
Consider $s_1=0$ and $y_1=1$. After dividing by $\kappa$ on both sides,  equation \eqref{eq:s1} becomes
\[(a-2b_0)(s_2^2(a+2b_0)+\kappa(a-2b_0)-4b_1\kappa s_2)=4t_2(-s_2+t_2).\]
Considering this as an equation in $t_2$ we see that this has solutions if \[s_2^2+(a-2b_0)(s_2^2(a+2b_0)+\kappa(a-2b_0)-4b_1\kappa s_2) \] is a square.
This is equivalent to \[f(s_2)=s_2^2(1+a^2-4b_0^2)+\kappa(a-2b_0)^2-4b_1\kappa s_2(a-2b_0)\] being a square for some $s_2\neq 0$. 
Note that if $a-2b_0=0$ this condition is trivially fulfilled. If $1+a^2-4b_0^2=1+(a+2b_0)(a-2b_0)=0$ then $a-2b_0\neq 0$ and in that case we see that $f(\frac{(1-\kappa)(a-2b_0)}{4b_1})=\kappa^2(a-2b_0)^2$ is a square, and $\frac{(1-\kappa)(a-2b_0)}{4b_1}\neq 0$.

If $1+a^2-4b_0^2\neq 0$ then $f(s_1)=f(s_2)$ if and only if $s_1+s_2=\frac{4b_1\kappa(a-2b_0)}{1+a^2-4b_0^2}$ so the image of $f$ has size $1+\frac{q-1}{2}$. Now $f(0)=\kappa(a-2b_0)^2$ is not a square since $-\kappa$ is not a square and $-1$ is a square. So since there are $(q-1)/2+1$ squares (including $0$) and $f(0)$ is not a square, we find that there is at least one $s_0\neq 0$ such that $f(s_0)$ is a square.




Now consider the case that $b_1=0$. We have $\norm(\beta)=b_0^2$, so $4\norm(\beta)-a^2=4b_0^2-a^2$ is a square
(because it lies in $\F_q$).
Consider $y_1=1$, then \eqref{eq:s1} becomes 
\[((a+2b_0)s_1^2+\kappa(a-2b_0))(s_2^2(a+2b_0)+\kappa(a-2b_0))=4t_2(s_1-s_2+t_2)\kappa.\]
Note that $\trace(\alpha+\beta)=a+2b_0\neq 0$ so we can multiply with $(a+2b_0)^2$ on both sides, and consider $\bar{s_1}=(a_2+b_0)s_1, \bar{s_2}=(a_2+b_0)s_2,\bar{t_2}=(a_2+b_0)t_2$, and look for solutions to the following equation:
\[(\bar{s_1}^2+\kappa(a^2-4b_0^2))(\bar{s_2}^2+\kappa(a^2-4b_0^2))=4\bar{t_2}(\bar{s_1}-\bar{s_2}+\bar{t_2})\kappa.\]
Let $\lambda=\kappa(a^2-4b_0^2)$ and recall that $\lambda$ is not a square.
As an equation in $t_2$, we have discriminant $4\kappa^2(\bar{s_1}-\bar{s_2})^2-4\kappa(\bar{s_1}^2+\lambda)(\bar{s_2}^2+\lambda)$ which needs to be a square. 
First assume that $\kappa\neq \lambda$. Let $\bar{s_2}=0$ then we need to find some $s\neq 0$ such that $g(s):=4\kappa(s^2(\kappa-\lambda)-\lambda^2)$ is a square.

First assume that $\kappa-\lambda$ is a non-square.
Consider $h(s)=g(s)/(4\kappa)$, we need to show that there is some $s\neq 0$ such that $h(s)$ is a non-square. Since $\kappa-\lambda$ is not a square, there is no $s$ with $h(s)=0$, and since $h(s)=h(t)\iff s=-t$, we find that $h$ takes on $1+(q-1)/2$ different nonzero values, so for at least one $s_0$, we have that $h(s_0)$ is a non-square, where $s_0\neq 0$ since $h(s_0)=-\lambda^2$ is a square. Now assume that $\kappa-\lambda=c^2$ is a nonzero square. Then $g(\lambda/c)=0$ which is a square.

Finally, assume that $\kappa-\lambda=0$, which happens if and only if $a^2-4b_0^2=1$. Let $\bar{s_1}=-\kappa$ and $\bar{s_2}=1$, then $\bar{s_1}\neq \bar{s_2}$ since $-\kappa$ is a non-square and $1$ is a square. The discriminant $4\kappa^2(\bar{s_1}-\bar{s_2})^2-4\kappa(\bar{s_1}^2+\lambda)(\bar{s_2}^2+\lambda)=4\kappa^2(\kappa+1)^2-4\kappa^2(\kappa+1)^2=0$ which is a square.
\end{proof}

We now can prove Theorem \ref{main2}, which we restate below:
\begin{quote}
    Let $q$ be an odd prime power, and let $\s$ be a set of $q^2+1$ pairwise noncollinear points of $\h(3,q^2)$ such that the following hold:
    \begin{enumerate}[(i)]
        \item There is a point $P\in\s$ and totally isotropic line $\ell$ through $P$, such that every point of $X$ incident
        with $\ell$ not equal to $P$ is collinear with a unique element of $\s\backslash\{P\}$. This yields
        a natural correspondence $F_\ell$ from $\ell$ to $\s$.
        \item Every Baer subline $b$ of $\ell$ containing $P$ maps under $F_\ell$ to $q+1$ 
        points of a plane of $\PG(3,q^2)$, incident with $P$.
        \item There is a point $Q$ such that for every $R\neq P,Q$ in $\s$, the triangle $PQR$ is in perspective.
        \item The plane spanned by three distinct points $R_1,R_2,R_3\in \s$ is nondegenerate.
    \end{enumerate}
    Then $\s$ is equivalent under the collineation group of $\h(3,q^2)$ to the Hermitian Veronesean $\mathcal{V}$.   
\end{quote}

\begin{proof}
By Lemma \ref{lemma:transitive}, we may suppose $P(0,0,0,1)$ and $Q(1,0,0,0)$.
By Lemma \ref{lemma:bijection}, its corollary, Lemma \ref{lemma:ct}, and Lemma \ref{lemma:flinear},
we can describe $\s$ as $S_{\alpha,\beta}$ for some $\alpha,\beta\in \F_{q^2}$, as
stated in Definition \ref{defn:salphabeta}. Suppose by way of contradiction that 
$\s$ is not equivalent to the Hermitian Veronesean $\mathcal{V}$ (under the collineation group of $\h(3,q^2)$).
Then by Lemma \ref{lemma:classical}, $\trace(\alpha)\ne 0$.
So, by Lemma \ref{lemma:tracenonzero}, there exist $t_1\in\F_{q^2}$, $t_2,t_3\in \F_q$, such that
$\trace[R_1,R_2,R_3]=0$ and $(t_1,t_2,t_3)$ is a distinct triple.
By Lemma \ref{lemma:degplane}, $R_1$, $R_2$, $R_3$
span a degenerate plane, a contradiction.
\end{proof}

Recall that a special set $\s$ of $\h(3,q^2)$ has the properties that a point not in $\s$ is collinear with $0$ or $2$ points of $\s$, and that every three points span a nondegenerate plane. Hence, for a special set, the conditions (i) and (iv) in Theorem \ref{main2} are satisfied and we find the following corollary:

\begin{corollary}\label{cor:sublinespecial}

Let $q$ be an odd prime power, and let $\s$ be a special set of $\h(3,q^2)$ such that the following hold:
    \begin{enumerate}[(i)]
        \item There is a point $P\in\s$ and a totally isotropic line $\ell$ through $P$ such that for every Baer subline $b$ of $\ell$ containing $P$, the points of $\s$ that are collinear with a point of $b$ are coplanar.
        \item There is a point $Q$ such that for every $R\neq P,Q$ in $\s$, the triangle $PQR$ is in perspective. 
    \end{enumerate}
    Then $\s$ is equivalent under the collineation group of $\h(3,q^2)$ to the Hermitian Veronesean $\mathcal{V}$. 
\end{corollary}

\section*{Acknowledgements} The first author thanks The University of Canterbury and
the Erskine Foundation for his `University of Canterbury Visiting Erskine Fellowship'
which supported this work.


\begin{thebibliography}{1}

\bibitem{Bamberg:2021aa}
J.~Bamberg, G.~Monzillo, and A.~Siciliano.
\newblock Pseudo-ovals of elliptic quadrics as {D}elsarte designs of
  association schemes.
\newblock {\em Linear Algebra Appl.}, 624:281--317, 2021.

\bibitem{Cossidente:2006aa}
A.~Cossidente, O.~H. King, and G.~Marino.
\newblock Special sets of the {H}ermitian surface and {S}egre invariants.
\newblock {\em European J. Combin.}, 27(5):629--634, 2006.

\bibitem{Grove}
L.~C. Grove.
\newblock {\em Classical groups and geometric algebra}, volume~39 of {\em
  Graduate Studies in Mathematics}.
\newblock American Mathematical Society, Providence, RI, 2002.

\bibitem{Hirschfeld:aa}
J.~W.~P. Hirschfeld and J.~A. Thas.
\newblock Arcs, caps and generalisations in a finite projective space.
\newblock arXiv:2503.06243.

\bibitem{kleidmanliebeck}
P.~Kleidman and M.~Liebeck.
\newblock {\em The Subgroup Structure of the Finite Classical Groups}, volume
  129 of {\em London Mathematical Society Lecture Note Series}.
\newblock Cambridge University Press, Cambridge, 1990.

\bibitem{Lavrauw:2023aa}
M.~Lavrauw, S.~Lia, and F.~Pavese.
\newblock On the geometry of the {H}ermitian {V}eronese curve and its
  quasi-{H}ermitian surfaces.
\newblock {\em Discrete Math.}, 346(10):Paper No. 113582, 20, 2023.

\bibitem{Segre65}
B.~Segre.
\newblock Forme e geometrie hermitiane, con particolare riguardo al caso
  finito.
\newblock {\em Ann. Mat. Pura Appl. (4)}, 70:1--202, 1965.

\bibitem{Shult:2005aa}
E.~E. Shult.
\newblock Problems by the wayside.
\newblock {\em Discrete Math.}, 294(1-2):175--201, 2005.

\bibitem{Thas:2011aa}
J.~A. Thas.
\newblock Generalized ovals in {${\rm PG}(3n-1,q)$}, with {$q$} odd.
\newblock {\em Pure Appl. Math. Q.}, 7(3, Special Issue: In honor of Jacques
  Tits):1007--1035, 2011.

\end{thebibliography}
\end{document}